\documentclass[a4paper,10pt]{article} 
\usepackage{amsmath, amsthm, amssymb}
\usepackage{url}

\usepackage{physics}
\usepackage{graphicx,lipsum}
\graphicspath{ {./downloads/} }

\usepackage[colorlinks,citecolor=red,urlcolor=blue,bookmarks=false,hypertexnames=true]{hyperref}
\usepackage{tikz}
\usetikzlibrary{calc}
\usetikzlibrary{shapes}
\usepackage[autostyle]{csquotes}
\makeatletter

\usepackage[colorlinks]{hyperref}
\usepackage[nameinlink,capitalize]{cleveref}
\newtheorem{theorem}{Theorem}[section]

\newtheoremstyle{named}{}{}{\itshape}{}{\bfseries}{.}{.5em}{\thmnote{#3's }#1}
\theoremstyle{named}

\theoremstyle{definition}
\newtheorem{definition}{Definition}[section]

\theoremstyle{remark}

\usepackage{xspace}
\usepackage[margin=1.0in]{geometry}

\usepackage{titlesec}

\usepackage{mathtools}

\DeclarePairedDelimiterX{\inp}[2]{\langle}{\rangle}{#1, #2}
\titleformat{\chapter}
  {\Large\bfseries} 
  {}                
  {0pt}            
  {\huge}

\title{Global Compatibility of Bi-Hamiltonian Structures on Three Dimensional Manifolds}

\author{Beste Madran \and Ender Abado\u{g}lu}

\begin {document}

\maketitle

\abstract { It is shown in \cite{yazar6} that a dynamical system defined by a nonvanishing vector field on an orientable three dimensional manifold is globally bi-Hamiltonian if and only if the first Chern class of the normal bundle of the given vector field vanishes, and the bi-Hamiltonian structure is globally compatible if and only if the Bott class of the complex codimension one foliation defined by the nonvanishing vector field vanishes. In this work, we constructed a dynamical system on $S^3$, which is globally bi-Hamiltonian, but the Hamiltonians are  not globally compatible.}


\section{INTRODUCTION}

In mathematics, as well as in physics and engineering applications, the concept of a three dimensional dynamical 
system has played central role \cite{k},\cite{a}. The dimension three has a special importance from the point of view of modeling real world phenomena, and also for some of its topological properties such as the existence of a unique differentiable structure \cite{16} and parallelizability \cite{17}. Constructing relations between dynamical
systems and geometric structures on a manifold would enable one to obtain some global information about the dynamical system by studying the global properties of the geometric structure. In this sense, the Hamiltonian structure seems to be one of the most powerful methods of investigation of the dynamical systems. The study of Hamiltonian structures is intensely developed over the last decades, especially after the monumental works of Magri \cite{13}, Weinstein \cite{15} and others. Basically, Hamiltonian structure relates a dynamical system on a manifold with one or more symplectic or Poisson structures on a manifold \cite{a}, \cite{19}. For an extensive review of multi-Hamiltonian theory of dynamical systems, we refer to \cite{18}. In this work, we will focus on bi-Hamiltonian structures defined by nonvanishing vector fields on orientable three dimensional manifolds. 

For general terminology and basic definitions about bi-Hamiltonian structures we refer to \cite{20}. In general, a pair of Poisson structures and their compatibility are fundamental features of bi-Hamiltonian systems. As it is well known, the existence and properties of compatible bi-Hamiltonian structure are intimately related with the complete integrability of the corresponding dynamical system. Also Bogoyavlenskii \cite{14} studied the Hamiltonian systems with incompatible Poisson bivectors (i.e., when $\left[ J_1 ,J_2 \right] _{SN}\neq0$ in even dimensions).

In  three dimensional case, it is shown in \cite{yazar5} that for a nonvanishing vector field on an orientable three dimensional manifold bi-Hamiltonian structure always
exists locally. However, it is shown in \cite{yazar6} that there are obstructions to the global existence and compatibility of bi-Hamiltonian structures. The tools that help to examine the global existence and global compatibility of Poisson structures are the Chern class of the normal bundle of the vector field, and the Bott class of the complex codimension one foliation defined by the vector field, respectively.

In this work, we try to answer the following question: Can we find an example of a dynamical system defined by a nonvanishing vector field on an orientable three dimensional manifold  which is globally bi-Hamiltonian but its bi-Hamiltonian structure is not globally compatible?
For this purpose, we first focus on the base manifold $M$. For the global existence of bi-Hamiltonian structure, the vanishing of the Chern class of the normal bundle of the vector field, which belongs to the second cohomology class of the base manifold, is required. Therefore our first condition on the base manifold will be $H^2(M)=0$. On the other hand, the obstruction to global compatibility is equivalent to vanishing of the Bott class of the complex codimension one
foliation defined by the given vector field, which belongs to the third cohomology class of the base manifold. Therefore, in order to construct an example of a bi-Hamiltonian system which is locally compatible but globally incompatible, one needs a base manifold with a vanishing second cohomology class but a nontrivial third cohomology class, i.e. we require $H^2(M)=0$ and $H^3(M) \neq 0$. Obviously, these requirements suggest that the base manifold could be cohomologous to $S^{3}$, and that's why we investigate the eigenforms of curl operator on $S^{3}$ to obtain an explicit example of global incompatibility of globally defined Hamiltonian functions.

In section 2, we will review the local and global existence theorems of bi-Hamiltonian structures. Then in section 3, we will construct our example on $S^{3}$ and give necessary proofs.

\section{THE EXISTENCE THEOREMS OF BI-HAMILTONIAN STRUCTURES}

In the analysis of dynamical systems, the quantities or properties that are invariant under the flow describing the dynamical system are important. Many properties of the system, such as Liouville integrability and stability are defined by or related with these invariants. The Hamiltonian function is a special invariant of the dynamical system on a manifold with a Poisson structure.

An autonomous dynamical system on a manifold $M$
\begin{equation}
\overset{\cdot }{x}\left( t\right) =v\left( x\left( t\right) \right) 
\end{equation}%
is generated by a vector field $v(x)$ on a manifold up to time reparametrization. 

A smooth function $I$ is called invariant under the flow of the vector field if
$\mathcal{L}_{v}I=0.$ If we relate the vector field to an invariant function via a Poisson structure $\mathbf{J}$  the vector field $v$ becomes a Poisson vector field and satisfy
\begin{equation}
\mathcal{L}_{v}\mathbf{J}=0.
\end{equation}

\begin{definition}
A Poisson structure on a manifold $M$ is the bilinear map $\left\{ \cdot
,\cdot \right\} $: $C^{\infty }\left( M\right) \times C^{\infty }\left(
M\right) \rightarrow C^{\infty }\left( M\right) $ satisfying

\begin{enumerate}
\item skew-symmetry condition, 
\begin{equation*}
\left\{ f,g\right\} =-\left\{ g,f\right\}   \label{e1}
\end{equation*}

\item Jacobi identity%
\begin{equation*}
\left\{ \left\{ f,g\right\} ,h\right\} +\left\{ \left\{ h,f\right\}
,g\right\} +\left\{ \left\{ g,h\right\} ,f\right\} =0 \label{e2}
\end{equation*}

\item Leibniz rule 
\begin{equation*}
\left\{ fg,h\right\} =f\left\{ g,h\right\} +g\left\{ f,h\right\}   \label{e3}.
\end{equation*}
\end{enumerate}
\end{definition}

By using the Jacobi identity and the pairing between $TM$  and $T^*M$ given by any metric on $M$,
\begin{equation}
X\left( f\right) =\left\langle df,X\right\rangle   \label{e4}
\end{equation}%
the Poisson bivector $\mathbf{J},$ the local section of $\Lambda ^{2}TM,$ is
defined by 
\begin{equation}
\left\{ f,g\right\} =\left\langle df\otimes dg,\mathbf{J}\right\rangle 
\label{e5}
\end{equation}%
and the Jacobi identity is given by 
\begin{equation}
\left[ \mathbf{J},\mathbf{J}\right] _{SN}=0  \label{e6}
\end{equation}%
where $\left[ \cdot ,\cdot \right] _{SN}$ is the Schouten-Nijenhuis bracket for bivector fields.

\begin{definition}

A dynamical system is said to be bi-Hamiltonian if there exist two different  Hamiltonian structures

\begin{equation}
v=\mathbf{J}_{1}\left( dH_{2}\right) =\mathbf{J}_{2}\left( dH_{1}\right) \label{e28}
\end{equation}
such that\textbf{\ }$\mathbf{J}_{1}$ and $\mathbf{J}_{2}$ are not multiples
of each other.

\end{definition}
 
The compatibility of the bi-Hamiltonian structure is defined by the compatibility of the Poisson structures of the bi-Hamiltonian system as follows: 
\begin{definition}
Two Poisson  structures $\mathbf{J}_{1}$ and $\mathbf{J}_{2}$ are called compatible if $\mathbf{J}%
_{1}(x)+c\mathbf{J}_{2}(x)$ is also a Poisson structure where c is an arbitrary constant \cite{20}.
\end{definition} 
This definition of compatibility is equivalent to the vanishing of the Schouten-Nijenhuis bracket of $\mathbf{J}%
_{1}$ and $\mathbf{J}_{2}$.
\subsection{The Local Existence Theorem of Bi-Hamiltonian Structures}

In this section, basic definitions related with a bi-Hamiltonian structure are given. For more information about Poisson structures of dynamical systems on three dimensional manifolds, we refer to \cite{yazar1}. 
Using any metric on the manifold, one can define a pairing between the tangent and cotangent bundles, and the bivector $\mathbf{J}$ can be considered as a map $%
\mathbf{J:}T^{\ast }M\rightarrow TM$.  By choosing a basis ${e_i}$ of $TM$, $\mathbf{J}$ could be represented by a 3x3 skew-symmetric 
tensor of type (2,0) which identifies an element of the Lie algebra $so\left( 3\right) $. Using
the isomorphism between $so\left( 3\right) $ and $
\mathbb{R}^{3}$ one could express the bivector $\mathbf{J}$ by a vector $J$, which is called the
Poisson vector

\begin{equation}
\mathbf{J=J}^{mn}e_{m}\wedge e_{n}=\varepsilon _{k}^{mn}{J}^{k}e_{m}\wedge
e_{n}
\end{equation}
and the vector field given by $J=J^{k}e_{k}$ is called the Poisson vector field on M.

In terms of the Poisson vector, Jacobi identity can be expressed by

\begin{equation}
\left( \nabla \times J\right) \cdotp J=0
\label{e25}
\end{equation}
and the bi-Hamiltonian structure given in $(\ref{e28})$
becomes
\begin{equation}
{v}=J_{1}\times\nabla H_{2}=J_{2}\times\nabla H_{1}.
\end{equation}

Since  $J_{1}$ and $J_{2}$ are nowhere multiples of each other by definition, they are linearly independent and we have 
\begin{equation}
J_{i}\cdot{v}=0.
\end{equation}

Now, defining
\begin{equation}
\widehat{e_{1}}=\frac{v}{\left\Vert v\right\Vert }
\end{equation}
we may extend it to a local orthonormal frame field $\left\{ \widehat{e_{1}},%
\widehat{e_{2}},\widehat{e_{3}}\right\} $. 
Orthonormality implies that $g(\widehat{e_{i}},\widehat{e_{j}})$ is constant and they satisfy
\begin{equation}
\left[ \widehat{e_{i}},\widehat{e_{j}}\right] =C_{ij}^{k}\left( x\right) 
\widehat{e_{k}}
\end{equation}
for some structure functions  $%
\left( C_{ij}^{k}\left( x\right) \right) $.

\begin{theorem}

\cite{yazar6} Any three dimensional dynamical system defined by a nonvanishing vector field on an orientable three dimensional manifold,

\begin{equation}
\overset{\cdot }{x}\left( t\right) =v\left( x\left( t\right) \right) 
\end{equation}
has a pair of compatible Poisson vectors

\begin{equation}
J_{i}=\alpha _{i}\left( \widehat{e_{2}}+\mu _{i}\widehat{e_{3}}\right) 
\end{equation}%
in which $\mu _{i}$'s are determined by the equation

\begin{equation}
\widehat{e_{1}}.\nabla \mu _{i}=-C_{31}^{2}-\mu _{i}\left(
C_{31}^{3}+C_{12}^{2}\right) -\mu _{i}^{2}C_{12}^{3} \label{e23}
\end{equation}%
and $\alpha _{i}$'s are determined by the equation

\begin{equation}
\widehat{e_{1}}.\nabla \ln \frac{\alpha _{i}}{\left\Vert v\right\Vert }%
=C_{31}^{3}+\mu _{i}C_{12}^{3}.
\end{equation}

Furthermore this dynamical system is a locally bi-Hamiltonian system with a pair of
local Hamiltonian functions determined by

\begin{equation}
J_{i}=\left( -1\right) ^{i+1}\phi \nabla H_{i}
\end{equation}%
where%
\begin{equation}
\phi =\frac{\alpha _{1}\alpha _{2}\left( \mu _{2}-\mu _{1}\right) }{%
\left\Vert v\right\Vert }.
\end{equation}
\end{theorem}

It follows from this theorem that, locally
it is always possible to find a pair of compatible Poisson structures such that the system defined by the nonvanishing vector field on a three dimensional manifold becomes bi-Hamiltonian.

\subsection{The Global Existence Theorem of Bi-Hamiltonian Structures}

In this section, the global properties of bi-Hamiltonian structures which are described in \cite{yazar6} will be summarized. In \cite{yazar6} it was shown that  a dynamical system defined by a nonvanishing vector field on an orientable three dimensional manifold is globally bi-Hamiltonian if and only if the first Chern class of the normal bundle of the given
vector field vanishes. Furthermore, the bi-Hamiltonian structure is globally compatible if and only if the Bott class of the complex codimension one
foliation defined by the given vector field vanishes.

We will associate the vector field with its normal bundle to investigate some global properties of the vector field $v$. Let E be the one dimensional subbundle of  $TM$ generated by $v$, and $Q = TM/E$ be the normal bundle of $v$. First, we are going to show that $Q$ is a complex line bundle over $M$.

\begin{definition}
A complex line bundle over a manifold $M$ is a
manifold $V$ and a smooth surjection $\pi :V\rightarrow M$ such that

\begin{enumerate}
\item Each fibre $\pi ^{-1}(m)=V_{m}$ is a  complex one dimensional  vector space.

\item Every $m\in M$ has an open neighborhood $U\in M$ for which there is a
diffeomorphism

\begin{equation}
\varphi :\pi ^{-1}(U)\rightarrow U\times 
\mathbb{C}
\end{equation}
such that $\varphi (V_{m})\subset \left\{ m\right\} \times \mathbb{C}$ \ for every $m$ and that moreover the map
\begin{equation}
\varphi |_{V_{m}}:V_{m}\rightarrow \left\{ m\right\} \times\mathbb{C}
\end{equation}
is a linear isomorphism.
\end{enumerate}
\end{definition}

Let ${\Lambda }:Q\rightarrow Q$ be the cross product with $\widehat{e_{1}}$, i.e. for any $v\in Q$
\begin{equation}
\Lambda(v)=\widehat{e_{1}}\times v.
\end{equation}
Then,
\begin{equation}
\Lambda^2(v)=\widehat{e_{1}}\times(\widehat{e_{1}}\times v)=-v
\end{equation}
and it defines a complex structure ${\Lambda }$ on
the fibers of $Q\rightarrow M$, and $Q$ becomes a complex line bundle over $M$.

\begin{definition}
A connection on $Q$ is an $\mathbb{R}$-bilinear map

$\ \ \ \ \ \ \ \ \ \ \ \ \ \ \ \ \ 
\begin{array}{ccc}
\ \ \ \ \nabla :\Gamma \left( TM\right) \times \Gamma \left( Q\right) \ \ \ 
& \rightarrow  & \Gamma \left( Q\right)  \\ 
\left( v,s\right)  & \mapsto  & \nabla _{v}s%
\end{array}$

satisfying
\begin{enumerate}

\item $\nabla _{v}\left( fs\right) =v\left( f\right) s+f\nabla _{v}s$

\item $\nabla _{fv}s=f\nabla _{v}s$.
\end{enumerate}
\end{definition} 

Now, let $\Omega$ be the volume form defined by the metric choosen, and we introduce the local  Poisson one-form $\mathbf{J}$ associated with a local Poisson bivector field $J$

\begin{equation}
\mathbf{J_i}=\iota _{J_i}\Omega.
\end{equation}%

The Jacobi identity for the Poisson form can then be written as

\begin{equation}
dJ_{i}\wedge J_{i}=0  \label{e7}
\end{equation}%
for each $i=1,2.$ Then, the compatibility condition becomes
\begin{equation}
J_1\wedge dJ_2+J_2\wedge dJ_1=0 \label{e24}
\end{equation}
and there exists a local
one-form $\gamma$, which defines a connection form on the normal bundle of
the given vector field,  such that  the Poisson one-forms satisfy
\begin{equation}
dJ_{i}=\gamma \wedge J_{i}  \label{e8}
\end{equation}%
for each $i=1,2$.

\subsubsection{Chern Class}

Chern classes are topological invariants associated with complex vector bundles on a smooth manifold. One can view Chern classes as (partial) obstructions to the vector bundle being trivial. In the case of a complex line bundle it is the only obstruction. Namely, a complex line bundle is trivial if and only if its Chern class vanishes, and we can find a global nonvanishing section of the complex line bundle.

Let $w=\iota _{v}\Omega$ be a two-form. It has constant rank 2 since $v$ is a nonvanishing vector field. Now, we try to find necessary and sufficient conditions under which

\begin{equation}
w=\iota _{v}\Omega =\phi dH_{1}\wedge dH_{2}
\end{equation}%
where the functions $\phi $, $H_{1}$ and $H_{2}$ are local, linearly independent and nowhere vanishing factors.
Such a decomposition may not exist globally. So, if $S_{w}$ is the subbundle of $TM$ on which $w$ is of maximal rank, then $S_{w}\cong Q$ defined above. Then,

\begin{equation}
w=p_{1}\wedge p_{2}
\end{equation}%
for which globally defined one-forms $p_{1}$ and $p_{2}$.

\begin{theorem}

\cite{d} Let $\Sigma$ be an $\mathbb R^n$ bundle over a connected base space $M$. Let $w$ be a two- form on $\Sigma$  of constant rank $2s$. 
Let $S_{w}$ be the subbundle of  $\Sigma$  on which $w$  is of maximal rank. $w$  decomposes if and only if
\begin{enumerate}
\item$S_{w}$ is a trivial bundle.
\item The representation of its normalization as a map $w_{1}:M\rightarrow  SO(2s)/U(s)$ arising from any trivialization of $S_{w}$  lifts to $SO(2s)$.
\end{enumerate}
\end{theorem}

As a consequence of this theorem, the criterion for the existence of a global bi-Hamiltonian structure is equivalent to the decomposition of the two-form $w$ into globally defined one-forms $p_{1}$ and $p_{2}$. This requires the vanishing of the first Chern class of $Q$ which means a complex line bundle $Q$ is trivial and hence it has a global section. 

Since Poisson one-forms and the integrability conditions are determined by the local solutions of $(\ref{e23})$, they are defined locally on each
local neighborhood. Let $\left\{ J_{i}^{p}\right\} $ and $\left\{ J_{i}^{q}\right\} $ be
the Poisson vector fields defined locally on $(U_{p},x_{p})$ and $(U_{q},x_{q})$
around points $p\in M$ and $q\in M$. 

If $c_{1}(Q)=0$, then two pairs of compatible
Poisson vector fields $\left\{ J_{i}^{p}\right\} $ and $\left\{
J_{i}^{q}\right\} $ on $U_{p}$ and $\ U_{q}$ respectively, are related on $%
U_{p}\cap U_{q}$ by%

\begin{equation}
\frac{J_{i}^{q}}{\left\Vert J_{i}^{q}\right\Vert }=\frac{J_{i}^{p}}{%
\left\Vert J_{i}^{p}\right\Vert }. \label{e22}
\end{equation}%

\begin{theorem}
\cite{yazar6} There exist two linearly independent global sections $%
\widehat{j_{i}}$ of $\ Q$ satisfying 

\begin{equation}
\widehat{j_{i}}\cdot\left( \nabla \times \widehat{j_{i}}\right) =0 \label{e21}
\end{equation}%
if and only if $c_{1}(Q)=0$, which are defined by $\widehat{j_{i}}$ $=\frac{J_{i}}{\left\Vert J_{i}\right\Vert }$.
\end{theorem}

In fact, the equation (\ref{e21}) is the equation defining the bi-Hamiltonian structure. Now, since all orientable three dimensional manifolds are parallelizable, $TM$ has a global frame. Since the vector field $v$ is also globally defined, it is possible to extend $v$ to a global frame. Then, the symmetry properties of (\ref{e21}) provides the Riccati equation given in (\ref{e23}), which is always equivalent to a second order linear O.D.E whose fundamental set always involves two independent functions, and these functions define two Poisson structures and hence the bi-Hamiltonian structure. Therefore, on each neighborhood we always have two independent solutions. On the intersection of two coordinate neighborhoods. Hence, there are four independent solutions of the Riccati equations on the intersection of coordinate neighborhoods. Now, by the symmetries of the Riccati equation, any four solutions of the Riccati equation are related to each other. Using this symmetry property of equations defining a bi-Hamiltonian structure, it is possible to obtain equation (\ref{e22}), which represents the symmetry of equations defining the bi-Hamiltonian structure. Under the condition of vanishing Chern class, this provides a way to construct the global solution. 

Indeed, this is a particular example of a more general method: Local geometric objects on manifolds  
are defined by certain equations called Lie equations whose local symmetry group defines the local Lie group leaving the geometric object invariant. We refer to \cite{ortacgil} for details about local Lie groups and their global properties. The solution of Lie equations on each coordinate neighborhood give the local geometric object on this neighborhood. When we want to extend this solution in a sheaf theoretic way, we obtained one solution from each neighborhood and expect them not to be equal but to coincide up to symmetry. In this way we not only extend the local geometric object itself to a global one but also we extend the local Lie group (which is nothing but the local symmetry group of the Lie equations defining the local geometric object) into a global Lie group. In our case, the symmetries of equation (\ref{e23}) or (\ref{e21}) allows us, under the condition of vanishing Chern class, to extend our solution to a global one using unit vectors in the direction of Poisson vectors.  
However, when we add the condition of compatibility given by (\ref{e24}) to this picture, this result will no longer be valid. Now the Lie equations for the locally compatible bi-Hamiltonian structure are given by (\ref{e21}) and (\ref{e24}) are different and therefore the local symmetry group the compatible bi-Hamiltonian structure is smaller. This restricts the solutions that could coincide and one may not extend the solutions, at least using the method proposed above. Indeed, we try to show that this is really the case and the existence of a globally compatible bi-Hamiltonian structure requires the vanishing of a secondary obstruction.     

\subsubsection{Bott Class}

Since $v$ is a nonvanishing vector field on $M$ then its integral curves are leaves of a one dimensional foliation. Since $Q=TM/E$ is a complex  line bundle on $M$, this foliation has complex codimension one. There are two secondary classes for complex codimension one foliations: the Godbillon–Vey class and the imaginary part of the Bott class \cite{yazar10}. Since in our assumption, the first Chern class of the complex normal bundle is trivial then we  compute the Bott class of this complex codimension one foliation.

By the Frobenius Theorem, the integrability of the distribution is equivalent to $\omega \wedge
d\omega =0$. That is, $d\omega =\gamma \wedge \omega $ for some one-form $\gamma$ and $\gamma$ is a connection on complex line bundle.

Assume that $Q$ and its dual $Q^{\ast }$ are trivial bundles. By $(\ref{e21})$, $Q^{\ast }$ has two global sections $\widehat{j_{i}}$ satisfying 

\begin{equation}
\widehat{j_{i}}\wedge d\widehat{j_{i}}=0\text{ and }\left\Vert \widehat{j_{i}%
}\right\Vert =1
\end{equation}%
and%

\begin{equation}
d\widehat{j_{i}}=\Gamma _{i}\wedge \widehat{j_{i}}
\end{equation}%
for globally defined $\Gamma _{i}$'s.

We state the following theorem provided in \cite{yazar6}, which will be essential for the construction of our example. 
\begin{theorem}\cite{yazar6}Let $\kappa $ be the curvature two-form of $Q$. There
exists a compatible pair of global Poisson structures if and only if%

\begin{equation}
\Xi =\left( \Gamma _{1}-\Gamma _{2}\right) \wedge \kappa
\end{equation}%
is exact. The cohomology class of  $\Xi $ vanishes if and only
if the Bott class of the complex codimension one foliation defined by the
nonvanishing vector field vanishes.
\end{theorem}
For the rest of this work, we try to construct a nonvanishing vector field on $S^3$ defining a nontrivial Bott class. For this purpose we will work on eigenforms of curl operator on $S^3$. Eigenforms of curl operator on $S^{3}$ provides the example of a locally compatible but globally incompatible pair of Hamiltonian structures  which are given by the vanishing of the cohomology classes 
\begin{equation}
\begin{array}{ccc}
c_{1} & = & \left[ d\gamma \right] \in H^{2}\left( S^{3}\right)  \\ 
&  &  \\ 
\beta  & = & \left[ \gamma \wedge d\gamma \right] \in H^{3}\left( S^{3}\right). 
\end{array}
\label{e9}
\end{equation}

Here the primary condition is the obstruction to the existence of 
bi-Hamiltonian structure and the secondary is the obstruction to global
compatibility. By the well known cohomology classes of  $S^{3}$, the first Chern class vanishes  since $H^{2}\left( S^{3}\right) =0$, which
implies that it is always possible to find a globally bi-Hamiltonian
structure but since  $H^{3}\left( S^{3}\right) =%
\mathbb{R}$ it may not always be possible to find a globally compatible pair.

\section{THE EXAMPLE ON $S^{3}$}

\subsection{The Eigenbasis of The Curl Operator on $S^{3}$}

We begin with the standard definition of the $S^{3}$ by using Cartesian coordinates 
\begin{equation}
S^{3}=\{(x^{1},x^{2},x^{3},x^{4})\in \mathbb{R}^{4}:(x^{1})^2+(x^{2})^2+(x^{3})^2+(x^{4})^2=R^2\}.
\end{equation}

Using Eulerian angles  $\tilde\varphi, \tilde\theta, \tilde\psi$, we can define the parametric representation of $S^{3}$  embedded in $\mathbb{R}^{4}$ as

\begin{equation}
\begin{array}{ccc}
x^{1} &=&R\cos (\nu\frac{\theta }{2})\cos(\nu \frac{\varphi +\psi }{2} )\\
&  &  \\ 
x^{2} &=&R\cos (\nu\frac{\theta }{2})\sin(\nu \frac{\varphi +\psi }{2}) \\
&  &  \\ 
x^{3} &=&R\sin(\nu \frac{\theta }{2})\cos (\nu\frac{\varphi -\psi }{2}) \\
&  &  \\ 
x^{4} &=&R\sin (\nu\frac{\theta }{2})\sin(\nu \frac{\varphi -\psi }{2})
\end{array}
\label{e26}
\end{equation}
where  $0\leq \nu\varphi <2\pi ,0\leq \nu\theta <\pi ,0\leq \nu\psi <4\pi $. The parameters $\varphi, \theta, \psi$ respectively represent the half-length of arcs which correspond to the Eulerian angles  $\tilde\varphi=\nu \varphi, \tilde\theta=\nu\theta, \tilde\psi=\nu\psi$ on the $S^{3}$ of radius $R=\frac2{\nu}$ \cite{yazar2}. For a more detailed  explanation of Eulerian angles, we refer to \cite{12}.

The metric on  $S^3$ of radius $R=2/\nu$ which is parameterized in terms of the (half) Eulerian arclength becomes
\begin{eqnarray}
ds^{2}=d^{2}\varphi +d^{2}\theta +d^{2}\psi +2\cos ( \nu\theta) d\varphi d\psi.
\end{eqnarray}

The special unitary group $SU(2)$ is diffeomorphic to $S^{3}$, and an element of $SU(2)$ is given by the matrix of $U$ 

\begin{equation}
U=\frac1{R}\left[ 
\begin{array}{c}
z_{1} \\ 
-\overline{z_{2}}%
\end{array}%
\begin{array}{c}
z_{2} \\ 
\overline{z_{1}}%
\end{array}
\right].
\end{equation}

The local description usually adopted is given by the conventional Eulerian arclength for a group element

\begin{equation}
U=exp(i{\nu\frac\varphi {2}}\sigma _{3})exp(i{\nu\frac\theta {2}}\sigma _{2})exp(i{\nu\frac\phi {2}}\sigma _{3})\end{equation}
where the Euler parameterization of an arbitrary matrix of $SU(2)$
transformation is defined by the appropriately chosen range for the Euler angles $\tilde\varphi , \tilde\theta, \tilde\psi$, and $\sigma _{i}$'s are the Pauli matrices which are the basis of $su(2)$. Using the normalization
\begin{equation}
\omega^i=-2tr(\tau_iU^{-1}dU)
\end{equation}
where $\tau_k=\frac{\sigma _{k}}{2i}$ and
by direct computation with the the Eulerian arclength representation $(\ref{e26})$ we reach the well known expressions for left-invariant one-forms which are also known as Maurer-Cartan forms. By using trigonometric half-angle formulas, we obtain our first form as
\begin{equation}
\omega ^{1} =\nu[\sin( \nu\theta) \cos (\nu\psi) d\varphi -\sin (\nu\psi) d\theta].
\end{equation}

For $\omega ^{2},\omega ^{3}$ all the calculation processes are almost the same. Therefore, we scale the unmodified coframe with the dimensionful factor $1/\nu$ \cite{yazar9} . This yields

\begin{eqnarray}
\omega ^{1} &=&\sin( \nu\theta) \cos (\nu\psi) d\varphi -\sin (\nu\psi) d\theta \label{e29}\\
\omega ^{2} &=&\sin( \nu\theta) \sin (\nu\psi) d\varphi +\cos (\nu\psi) d\theta \label{e30}\\
\omega ^{3} &=&\cos ( \nu\theta) d\varphi +d\psi \label{e31}. 
\end{eqnarray}

The determinant of this metric is
\begin{equation}
\det g_{\mu \upsilon }=\left( \det g^{\mu \upsilon }\right) ^{-1}=\left(
\det \omega _{\mu }^{i}\right) ^{2}=\sin ^{2}( \nu\theta) =\left\vert g\right\vert.
\end{equation}
The invariant volume element on $SU(2)$ \cite{yazar7} ,
\begin{equation}
d\mu(g)=\delta(R-1)dx^1dx^2dx^3dx^4
\end{equation}
calculated with the differentials of $x^i$'s takes the form,
\begin{equation}
d\mu(g)=\frac1{16\pi^2}sin(\nu\theta)d\varphi d\theta d\psi
\end{equation}
which is normalized so that

\begin{equation}
\int_{SU\left( 2\right)}d\mu(g)=1
\end{equation}
The invariant volume element on $SU(2)$ is related with the volume element of $S^3$, which we denote by $dVol$ by
\begin{equation}
dVol=2\pi^2d\mu(g).
\end{equation}
The invariant integral on $S^{3}$ expressed in terms of Eulerian arclength is

\begin{equation}
\int_{SU\left( 2\right)}f(g)d\mu(g)=\frac1{16\pi^2}\int_{0}^{\frac{2\pi}\nu}\int_{0}^{\frac{\pi}\nu }\int_{0}^{\frac{4\pi}\nu}f(\varphi,\theta,\psi)\sin( \nu\theta) d\varphi d\theta d\psi .
\end{equation}

Furthermore according to the well known 
the Maurer-Cartan equation in orthonormal basis
\begin{equation}
d\omega ^{k}+\frac{1}{2}\epsilon _{ijk}\omega ^{i}\wedge \omega ^{j}=0
\end{equation}
for $i,j,k=1,2,3$ gives 
\begin{eqnarray}
d\omega ^{1} &=&-\nu\omega ^{2}\wedge \omega ^{3} \\
d\omega ^{2} &=&-\nu\omega ^{3}\wedge \omega ^{1} \\
d\omega ^{3} &=&-\nu\omega ^{1}\wedge \omega ^{2}.
\end{eqnarray}%

The coframe $\{ \omega^1,\omega^2,\omega^3 \}$ determines a unique orientation and satisfying the
Hodge duality relations such as 
\begin{eqnarray}
\ast \omega ^{1} &=&\omega ^{2}\wedge \omega ^{3} \\
\ast \omega ^{2} &=&\omega ^{3}\wedge \omega ^{1} \\
\ast \omega ^{3} &=&\omega ^{1}\wedge \omega ^{2}.
\end{eqnarray}

And using inverse matrix  $\omega^{i},i=1,2,3$ the coordinate frame looks like
\begin{eqnarray}
d\varphi &=&\frac{\sin (\nu\psi) }{\sin ( \nu\theta) }\omega ^{1}+\frac{\cos (\nu\psi) }{%
\sin ( \nu\theta) }\omega ^{2} \\
d\theta &=&\cos (\nu\psi)\omega ^{1}-\sin (\nu\psi) \omega ^{2} \\
d\psi &=&-\sin (\nu\psi) \cot ( \nu\theta) \omega ^{1}-\cos (\nu\psi)\cot ( \nu\theta)\omega
^{2}+\omega ^{3}.
\end{eqnarray}

The curl operator on one-forms is $\ast d : \Omega^1(M) \rightarrow\Omega^1(M)$, where $\ast$ is the Hodge star operator and $d$ is the exterior derivative. The eigenform of the curl is any one-form $\omega $ satisfying
\begin{equation}
\ast d\omega =\lambda \omega   \label{e17}
\end{equation}%
for some $\lambda\in\mathbb{R}$.

On $S^3$, given a nonvanishing curl eigenform $\omega ^{3}$ with nonzero eigenvalue $\lambda$, since the volume form on a compact manifold cannot be exact unless $\lambda
=0$, it easily follows that

\begin{equation}
\begin{array}{ccc}
\ast d\omega ^{3} &=&\lambda \omega ^{3} \\
&  &  \\ 
\ast \left(-\nu \omega ^{1}\wedge \omega ^{2}\right) &=&\lambda \omega ^{3} \\
&  &  \\ 
-\nu\omega ^{3} &=&\lambda \omega ^{3}
\end{array}
\label{e18}
\end{equation}
which implies 
\begin{equation}
\lambda =-\nu.
\end{equation}

Now, we will investigate the bi-Hamiltonian structure of the dynamical system defined by $\omega ^{3}$. Let us briefly recall that locally, every pair of Poisson vectors can be made compatible.
Hence there exists a local one-form $\gamma$, which is also a connection form on the normal bundle of the given vector field, such that  the Poisson one-forms satisfy
\begin{equation}
dJ_{i}=\gamma \wedge J_{i}  \label{e19}
\end{equation}
for each $i = 1,2$.

One could easily find that the two local Hamiltonian functions are
\begin{eqnarray}
H_{1}=\nu\frac{\varphi +\psi }{2}\\
H_{2}=\nu\frac{\varphi -\psi }{2}
\end{eqnarray}
and our local Poisson structures are
\begin{eqnarray}
J_{1} &=&d\left( \nu\frac{\varphi +\psi }{2}\right) \\
J_{2} &=&d\left( \nu\frac{\varphi -\psi }{2}\right) 
\end{eqnarray}
which satisfy
\begin{equation}
dJ_{i}=\gamma \wedge J_{i}  \label{e20}
\end{equation}
where $\gamma =0$.

In order to write local Poisson structures by using orthonormal bases $(\ref{e29})-(\ref{e31})$ gives
\begin{eqnarray}
J_{1} &=&\frac{\nu}{2}\left[ \frac{\sin (\nu\psi) }{\sin ( \nu\theta) }\left( 1-\cos
( \nu\theta) \right) \omega ^{1}+\frac{\cos(\nu \psi) }{\sin ( \nu\theta) }\left( 1-\cos
( \nu\theta) \right) \omega ^{2}+\omega ^{3}\right] \\
J_{2} &=&\frac{\nu}{2}\left[ \frac{\sin (\nu\psi) }{\sin ( \nu\theta) }\left( 1+\cos
( \nu\theta) \right) \omega ^{1}+\frac{\cos (\nu\psi) }{\sin ( \nu\theta) }\left( 1+\cos
( \nu\theta) \right) \omega ^{2}-\omega ^{3}\right].
\end{eqnarray}

Recall that we are working in a chart 
$0\leq \nu\varphi <2\pi ,0\leq \nu\theta <\pi ,0\leq \nu\psi <4\pi $. But in this chart 
\begin{equation}
\lim_{\nu\theta \to \pi} J_{1}\to\infty,
\lim_{\nu\theta \to 0} J_{2}\to\infty.
\end{equation} 

\subsection{The Obstruction to Global Compatibility}

The unit Poisson vector fields are
\begin{equation}
\widehat{j_{i}}=\frac{J_{i}}{\left\Vert J_{i}\right\Vert }
\end{equation} 
where

\begin{eqnarray}
\left\Vert J_{1}\right\Vert &=&\frac{\nu}{2}\sec (\nu\frac{\theta }{2})\\
\left\Vert J_{2}\right\Vert& =&\frac{\nu}{2}\csc (\nu\frac{\theta }{2})
\end{eqnarray}
and therefore the global Poisson structures are
\begin{eqnarray}
\widehat{j_{1}} &=&\sin (\nu\psi) \sin(\nu \frac{\theta }{2})\omega ^{1}+\cos (\nu\psi)
\sin(\nu \frac{\theta }{2})\omega ^{2}+\cos (\nu\frac{\theta }{2})\omega ^{3} \\
\widehat{j_{2}} &=&\sin(\nu \psi) \cos(\nu \frac{\theta }{2})\omega ^{1}+\cos(\nu \psi)
\cos (\nu\frac{\theta }{2})\omega ^{2}-\sin(\nu \frac{\theta }{2})\omega ^{3}
\end{eqnarray}
with differentials 

\begin{equation}
d\widehat{j_{1}}=\frac{\nu}{2}\left[ \sin (\nu\psi) \sin(\nu \frac{\theta }{2})\omega
^{2}\wedge \omega ^{3}+\cos(\nu \psi) \sin(\nu \frac{\theta }{2})\omega ^{3}\wedge
\omega ^{1}-\frac{\sin ^{2}(\nu\frac{\theta }{2})}{\cos (\nu\frac{\theta }{2})}\omega
^{1}\wedge \omega ^{2}\right]
\end{equation}%

\begin{equation}
d\widehat{j_{2}}=\frac{\nu}{2}\left[ \sin (\nu\psi) \cos(\nu \frac{\theta }{2})\omega
^{2}\wedge \omega ^{3}+\cos(\nu \psi) \cos(\nu \frac{\theta }{2})\omega ^{3}\wedge
\omega ^{1}-\frac{\cos ^{2}(\nu\frac{\theta }{2})}{\sin (\nu\frac{\theta }{2})}\omega
^{1}\wedge \omega ^{2}\right].
\end{equation}.

Now, the Jacobi condition  
\begin{equation}
d\widehat{j_{i}}\wedge \widehat{j_{i}}=0
\end{equation}
is satisfied, and this implies that there exist two global connections 

\begin{equation}
d\widehat{j_{i}}=\Gamma _{i}\wedge \widehat{j_{i}} \label{e27}
\end{equation}
so we need to calculate connections $\Gamma _{i}$ for $i=1,2$.

Let
\begin{equation}
\Gamma _{i}=a_{i1}\omega ^{1}+a_{i2}\omega ^{2}+a_{i3}\omega ^{3}
\end{equation}
be any connection of one-form where $a_{ij}$'s are arbitrary. It is easy to check that the system of equations $(\ref{e27})$ has infinitely many
solutions. To find one, suppose $a_{13}=t_{1}$ then,
 
\begin{eqnarray}
a_{11} &=&\tan (\nu\frac{\theta }{2})\left( -\frac{\nu}{2}\cos(\nu \psi) +t_{1}\sin (\nu\psi)\right)  \\
a_{12} &=&\tan (\nu\frac{\theta }{2})(\frac{\nu}{2}\sin (\nu\psi )+t_{1}\cos (\nu\psi )).
\end{eqnarray}%
So our first connection becomes

\begin{equation}
\Gamma _{1}=\tan (\nu\frac{\theta }{2})\left[ \left( -\frac{\nu}{2}\cos(\nu \psi) 
+t_{1}\sin (\nu\psi)\right) \omega ^{1}+(\frac{\nu}{2}\sin (\nu\psi )+t_{1}\cos (\nu\psi )
)\omega ^{2}\right] +t_{1}\omega ^{3}.
\end{equation}
 
For $\Gamma _{2},$ all calculations and results are similar, and we have
\begin{equation}
\Gamma _{2}=\cot (\nu\frac{\theta }{2})\left[ \left( \frac{\nu}{2}\cos(\nu \psi) 
-t_{2}\sin (\nu\psi)\right) \omega ^{1}+(-\frac{\nu}{2}\sin (\nu\psi )-t_{2}\cos (\nu\psi )
)\omega ^{2}\right] +t_{2}\omega ^{3}.
\end{equation}

Now, we will to calculate  the following three-form 
\begin{equation}
\Xi =\left( \Gamma _{1}-\Gamma _{2}\right) \wedge \widehat{j_{1}}\wedge 
\widehat{j_{2}}.
\end{equation}
We begin with subtracting the $\Gamma _{i}$ 's from each other   
\begin{eqnarray*}
\left( \Gamma _{1}-\Gamma _{2}\right) 
&=&\left[ -\frac{\nu}{\sin (\nu\theta) }\cos(\nu \psi) +\left( t_{1}\tan(\nu \frac{\theta }{2%
})+t_{2}\cot(\nu \frac{\theta }{2})\right) \sin (\nu\psi) \right] \omega ^{1}
\\
&&+\left[ \frac{\nu}{\sin (\nu\theta) }\sin(\nu \psi) +\left( t_{1}\tan(\nu \frac{\theta }{2%
})+t_{2}\cot(\nu \frac{\theta }{2})\right) \cos (\nu\psi) \right] \omega ^{2}\\
&&+\left[t_{1}+t_{2}
\right] \omega ^{3}
\end{eqnarray*}
and calculating the wedge product of $\widehat{j_{i}}$'s
\begin{equation}
\widehat{j_{1}}\wedge \widehat{j_{2}}=-\cos (\nu\psi) \omega ^{2}\wedge \omega
^{3}+\sin(\nu \psi) \omega ^{3}\wedge \omega ^{1}
\end{equation}
and multiplying both sides by $\left( \Gamma _{1}-\Gamma _{2}\right) $ we get
\begin{equation}
\Xi =\left( \Gamma _{1}-\Gamma _{2}\right) \wedge \widehat{j_{1}}\wedge \widehat{%
j_{2}}=\frac{\nu}{\sin(\nu \theta) }\omega ^{1}\wedge \omega ^{2}\wedge \omega ^{3}.
\end{equation}
\begin{theorem} 
The three-form $\Xi $ is not exact.
\end{theorem}
\begin{proof}
In \cite{yazar7} for the compact group $SU(2)$ there exists left and right invariant measure $d\mu(g)$ satisfying
\begin{equation}
\int_{SU(2)}f(g'g)d\mu(g)=\int_{SU(2)}f(gg')d\mu(g)=\int_{SU(2)}f(g)d\mu(g)
\end{equation}
where $f(g)$ is a continuous function on $SU(2)$. Then normalized Haar integral of the function $f$ can be simplified as
\begin{equation}
\int_{SU(2)}f(g)d\mu(g)=\int_{0}^{\frac{\pi}\nu }f(\theta)sin(\nu\theta) d\theta.
\end{equation}
Recall that we are working in a chart 
$0\leq \nu\varphi <2\pi ,0\leq \nu\theta <\pi ,0\leq \nu\psi <4\pi. $

\begin{eqnarray*}
\int _{S^3} \frac{\nu}{\sin(\nu \theta) }\omega ^{1}\wedge \omega ^{2}\wedge \omega ^{3}
&=&\int_{0}^{\frac{2\pi}\nu}\int_{0}^{\frac{\pi}\nu }\int_{0}^{\frac{4\pi}\nu}\frac{\nu}{\sin(\nu \theta) }\sin(\nu\theta) d\varphi \wedge d\theta \wedge d\psi \\
&=&\nu\int_{0}^{\frac{2\pi}\nu}\int_{0}^{\frac{\pi}\nu }\int_{0}^{\frac{4\pi}\nu}d\varphi \wedge d\theta \wedge d\psi \\
&=&\frac{8\pi ^{3}}{\nu^{2}}.
\end{eqnarray*}
\end{proof} 

This proves the global incompatibility of the bi-Hamiltonian structure. However, for a double check, we are also going to show that the Bott class of the complex codimension one foliation does not vanish. The general form of connection can be written as 

\begin{equation}
\Gamma =\left[ 
\begin{array}{cc}
\Gamma _{1}+b_{11}\widehat{j_{1}} & b_{12}\widehat{j_{2}} \\ 
b_{21}\widehat{j_{1}} & \Gamma _{2}+b_{22}\widehat{j_{2}}%
\end{array}%
\right]
\end{equation}
where $b_{ij}$'s are arbitrary for each $i,j=1,2$.

Then, its curvature becomes
\begin{equation}
\kappa =\left[ 
\begin{array}{cc}
d\Gamma _{1}+(db_{11}+b_{11}\Gamma _{1}+b_{12}b_{21}\widehat{j_{2}})\wedge 
\widehat{j_{1}} & (db_{12}+b_{12}\Gamma _{1}+b_{11}b_{12}\widehat{j_{1}}%
)\wedge \widehat{j_{2}} \\ 
(db_{21}+b_{21}\Gamma _{2}+b_{21}b_{22}\widehat{j_{2}})\wedge \widehat{j_{1}}
& d\Gamma _{2}+(db_{22}+b_{22}\Gamma _{2}+b_{21}b_{12}\widehat{j_{1}})\wedge 
\widehat{j_{2}}%
\end{array}%
\right].
\end{equation}

By a simple calculation, we get the trace of curvature 
\begin{eqnarray}
tr\left( \kappa \right) =d\left( \Gamma _{1}+\Gamma _{2}+b_{11}\widehat{j_{1}}%
+b_{22}\widehat{j_{2}}\right)
\end{eqnarray}
and the Chern class is trivial
\begin{eqnarray}
c_{1} &=&\left[d\left( \Gamma _{1}+\Gamma _{2}+b_{11}\widehat{j_{1}}+b_{22}%
\widehat{j_{2}}\right)\right]=0
\end{eqnarray}
and the Bott class becomes%
\begin{eqnarray*}
\left[ \beta \right]  &=& \left[h_{1}\wedge c_{1} \right]\\
&=&\left[\left( \Gamma _{1}+\Gamma _{2}\right) \wedge d\left( \Gamma _{1}+\Gamma
_{2}\right) +\left( \Gamma _{1}+\Gamma _{2}\right) \wedge d\left( b_{11}%
\widehat{j_{1}}+b_{22}\widehat{j_{2}}\right) \right]  \\
&+&\left[\left( b_{11}\widehat{j_{1}}+b_{22}\widehat{j_{2}}\right) \wedge d\left(
\Gamma _{1}+\Gamma _{2}\right) +\left( b_{11}\widehat{j_{1}}+b_{22}\widehat{%
j_{2}}\right) \wedge d\left( b_{11}\widehat{j_{1}}+b_{22}\widehat{j_{2}}%
\right) \right]
\end{eqnarray*}%
where $dh_{1}=c_{1}$.

Now, choose $b_{11},b_{22}=1$ then try to calculate  Bott class. It consist of four terms,
we will try to calculate in order.

Firstly we added the $\Gamma _{i}$  's from each other 

\begin{eqnarray*}
\Gamma _{1}+\Gamma _{2} &=&\left[ \nu\cot(\nu \theta) \cos (\nu\psi) +\left( t_{1}\tan(\nu \frac{\theta }{2})%
-t_{2}\cot (\nu\frac{\theta }{2})\right) \sin (\nu\psi) \right] \omega ^{1} \\
&+&\left[ -\nu\cot (\nu\theta) \sin(\nu \psi) +\left( t_{1}\tan (\nu\frac{\theta }{2})%
-t_{2}\cot (\nu\frac{\theta }{2})\right) \cos (\nu\psi )\right] \omega ^{2}\\
&+&\left(t_{1}+t_{2}\right) \omega ^{3}
\end{eqnarray*}
then take the derivative of $\Gamma _{1}+\Gamma _{2}$, it is too long calculation but simple so we obtain  

\begin{eqnarray}
d\left( \Gamma _{1}+\Gamma _{2}\right)  &=&0.
\end{eqnarray}

Hence
\begin{equation}
\left( \Gamma _{1}+\Gamma _{2}\right) \wedge d\left( \Gamma _{1}+\Gamma
_{2}\right) =0.
\end{equation}

The second term of Bott class is $\left( \Gamma _{1}+\Gamma _{2}\right) \wedge d\left( \widehat{j_{1}}+\widehat{j_{2}}\right)$. Firstly we calculate the derivative of the sum of Poisson structures
\begin{eqnarray*}
 d\left( \widehat{j_{1}}+\widehat{j_{2}}\right)&=& \frac{\nu}{2}\left[ \sin(\nu \psi) \left( \sin(\nu \frac{\theta }{2})+\cos(\nu \frac{\theta }{2})\right)\right] \omega^{2}\wedge \omega ^{3} \\ 
&+& \frac{\nu}{2}\left[ \cos(\nu \psi) \left( \sin(\nu \frac{\theta }{2})+\cos (\nu\frac{\theta }{2})\right)\right]
\omega ^{3}\wedge \omega ^{1}\\
&+& \frac{\nu}{2}\left[ \left( -\frac{\sin ^{2}(\nu\frac{\theta }{2})}{\cos 
(\nu\frac{\theta }{2})}+\frac{\cos ^{2}(\nu\frac{\theta }{2})}{\sin (\nu\frac{\theta }{2})}
\right)\right] \omega ^{1}\wedge \omega ^{2}.
\end{eqnarray*}

Hence 
\begin{equation}
\left( \Gamma _{1}+\Gamma _{2}\right) \wedge d\left( \widehat{j_{1}}+%
\widehat{j_{2}}\right) 
=\frac{\nu}{2}\left[ t_{1}\frac{1}{\sin(\nu \frac{\theta }{2})}-t_{2}\frac{1}{%
\cos(\nu \frac{\theta }{2})}\right] \omega ^{1}\wedge \omega ^{2}\wedge \omega
^{3} \\.
\end{equation}

The third term of Bott class becomes 
\begin{equation}
\left( \widehat{j_{1}}+\widehat{j_{2}}\right) \wedge d\left( \Gamma
_{1}+\Gamma _{2}\right) =0.
\end{equation}

Moreover the last term is

\begin{equation}
\left( \widehat{j_{1}}+\widehat{j_{2}}\right) \wedge \left( d\widehat{j_{1}}%
+d\widehat{j_{2}}\right) =\frac{\nu}{\sin(\nu \theta) }\omega ^{1}\wedge \omega
^{2}\wedge \omega ^{3}.
\end{equation}

Consequently, the Bott class is equal to
\begin{equation}
\left[\beta \right] =\left[ \left[\frac{\nu}{2}\left[ t_{1}\frac{1}{\sin(\nu \frac{\theta }{2})}-t_{2}\frac{1}{%
\cos(\nu \frac{\theta }{2})}\right] +\frac{\nu}{\sin(\nu \theta) }\right]\omega ^{1}\wedge \omega
^{2}\wedge \omega ^{3} \right].
\end{equation}

Therefore   $t_{i}$'s are arbitrary so call these two terms $t_{1}=\sin(\nu \frac{\theta }{2})$, $t_{2}=\frac{1}{\sin 
(\nu\frac{\theta }{2})}$ 
and showed that 
\begin{equation}
\beta =\frac{\nu}{2}\omega ^{1}\wedge \omega
^{2}\wedge \omega ^{3}.
\end{equation}

Calculating the integral of $\beta $

\begin{eqnarray}
\int _{S^3} \frac{\nu}{2}\omega ^{1}\wedge \omega ^{2}\wedge \omega ^{3}&=&
 \frac{\nu}{2}\int_{0}^{\frac{2\pi}\nu}\int_{0}^{\frac{\pi}\nu }\int_{0}^{\frac{4\pi}\nu}\sin (\nu\theta)
d\varphi d\theta d\psi \\
&=&\frac{8\pi ^{2}}{\nu^{2}}
\end{eqnarray}
which implies that the Bott class is not trivial either.

\end{document}